\documentclass[imammb]{ima-authoring-template}%
\makeatletter
\newcommand{\@copyrightyear}{202x}
\makeatother

\graphicspath{{Fig/}}

\usepackage{graphicx} % Required for inserting images
\usepackage{amsmath}
\usepackage{amsfonts}
\usepackage{amssymb}
\usepackage{amsthm}
\usepackage{hyperref}
\usepackage{caption}
\usepackage{float}
\usepackage[normalem]{ulem}
\usepackage{cancel}
\usepackage{booktabs}
\usepackage{multirow}
\usepackage{multicol}

\usepackage{natbib} % for author-year citations

\newcommand{\train}[1]{\ensuremath{\mathcal{#1}_\text{tr}}}
\newcommand{\test}[1]{\ensuremath{\mathcal{#1}_\text{te}}}

\newcommand{\ensuretext}[1]{\ifmmode\text{#1}\else#1\fi}
\newcommand{\qty}[2]{#1\,#2}
\newcommand{\qtyrange}[3]{#1--#2\,#3}
\newcommand{\si}[1]{#1}

\newcommand{\Pa}{\ensuretext{Pa}}
\newcommand{\second}{\ensuretext{s}}
\newcommand{\N}{\ensuretext{N}}
\newcommand{\meter}{\ensuretext{m}}
\newcommand{\kg}{\ensuretext{kg}}
\newcommand{\micro}{\ensuremath{\mu}}
\newcommand{\milli}{\ensuretext{m}}
\newcommand{\minute}{\ensuretext{min}}
\newcommand{\mol}{\ensuretext{mol}}

\newcommand{\um}{\micro\meter}
\newcommand{\Osm}{\ensuretext{Osm}}

\newcommand{\per}{\ensuremath{\,\!/\!\,}}
\newcommand{\tothe}[1]{\ensuremath{^{#1}}}
\newcommand{\squared}{\tothe{2}}
\newcommand{\cubed}{\tothe{3}}

\numberwithin{equation}{section}

\begin{document}

\title{Operator learning for models of tear film breakup}

\author{Qinying Chen\footnote{These authors contributed equally to this work.},
Arnab Roy\footnotemark[1], and Tobin A. Driscoll
\address{\orgdiv{Department of Mathematical Sciences},
\orgname{University of Delaware},
\orgaddress{\street{Newark}, \postcode{19716}, \state{DE}, \country{USA}}}}

\authormark{Chen et al.}

\corresp[*]{Corresponding author: Arnab Roy, \href{mailto:arnabroy@udel.edu}{arnabroy@udel.edu}}

\received{Date}{0}{Year}
\revised{Date}{0}{Year}
\accepted{Date}{0}{Year}

\abstract{Tear film (TF) breakup is a key driver of understanding dry eye disease, yet estimating TF thickness and osmolarity from fluorescence (FL) imaging typically requires solving computationally expensive inverse problems. We propose an operator learning framework that replaces traditional inverse solvers with neural operators trained on simulated TF dynamics. This approach offers a scalable path toward rapid, data-driven analysis of tear film dynamics.}


\maketitle

\section{Introduction}

A thin layer of liquid known as the tear film (TF) lies across the surface of the cornea of a human eye~\citep{Doane80}. This dynamic film provides lubrication, protects against microbial invasion, maintains a smooth optical interface, and delivers essential oxygen and nutrients to the cornea \citep{willcoxTFOSDEWSII2017}. The TF has three layers: a lipid layer that is 20 to 100 nm thick \citep{braunDynamicsFunctionTear2015}, an aqueous layer consisting mainly of water \citep{hollyFormationRuptureTear1973} that is a few microns thick, and a half-micron thick mucin layer called the glycocalyx that is on the ocular surface \citep{king-smithThicknessTearFilm2004}. The lipid layer slows evaporation of water from the TF \citep{mishimaOilyLayerTear1961}, and a healthy glycocalyx facilitates the fluid movements of the ocular surface \citep{gipsonDistributionMucinsOcular2004}. The lacrimal gland supplies the majority of the water from the aqueous layer near the temporal canthus \citep{darttNeuralRegulationLacrimal2009}. Osmosis supplies water from the ocular epithelia \citep{braunDynamicsTear2012}.

When the tear film becomes unstable and breaks up, a process known as tear breakup (TBU), it can leave patches of the ocular surface exposed, triggering irritation, inflammation, and visual disturbances \cite{king2018mechanisms}. Persistent or frequent tear breakup is recognized as a major factor in the onset and progression of dry eye disease (DED), a multifactorial condition affecting millions of individuals worldwide that compromises both the comfort and visual function of the eye~\citep{lempDefinitionClassificationDry2007, nelsonTFOSDEWSII2017, stapletonDEWSIIepi2017}.
One common subtype, evaporative dry eye (EDE), is primarily driven by excessive water loss from the tear film, often due to dysfunction in the lipid layer that normally retards evaporation \citep{linDryEyeDisease2014, OCEANreport2013}. TBU happens when a dry spot appears on the eye \citep{nornMICROPUNCTATEFLUORESCEINVITAL1970} and is often evaporation-driven \citep{lempDefinitionClassificationDry2007,willcoxTFOSDEWSII2017}. 

Fluorescein (FL) imaging is one of the most widely used techniques for studying tear-film dynamics and tear-film breakup. Following instillation of sodium fluorescein, the aqueous layer of the tear film emits fluorescence whose intensity depends on both tear-film thickness and fluorescein concentration \citep{webberContinuousFluorophotometricMethod1986,nicholsThinningRatePrecorneal2005a}. FL imaging has been used to measure tear-film breakup times, visualize the formation and progression of breakup regions, and investigate tear-film thinning dynamics \citep{nornMICROPUNCTATEFLUORESCEINVITAL1970,begleyQuantitativeAnalysisTear2013,liTearFilmDynamics2014}. The interpretation of FL images, however, is not straightforward, because the observed intensity depends on both tear-film thickness and fluorescein concentration. In particular, the relationship between fluorescence intensity and tear-film thickness varies between the dilute and self-quenching concentration regimes \citep{webberContinuousFluorophotometricMethod1986, King-SmithIOVS13a,braunModelTearFilm2014}. Previous modeling studies incorporated fluorescein transport and intensity calculations to connect FL image data with underlying tear-film dynamics. Using combined experimental and modeling approaches, these studies showed that evaporation, osmolarity changes, and tangential flow can all influence the appearance of breakup regions in FL images and that advection of fluorescein may complicate direct interpretation of intensity measurements \citep{braunDynamicsFunctionTear2015, braunTearFilm2018, lukeParameterEstimation2020}.

Measurement of TF dynamics in situ is challenging. Evaporation has been identified as a primary driver of the dynamics, with models incorporating spatially varying lipid layers, osmolarity transport in the aqueous layer, and osmosis across the cornea \citep{PengEtal2014,braunDynamicsFunctionTear2015,braunTearFilm2018}. These studies showed that osmolarity diffusion prevents osmosis from halting thinning, a phenomenon absent in spatially uniform models \cite{braunDynamicsTear2012}. Simpler models later included fluorescein transport, enabling direct comparison with imaging experiments, and one-dimensional partial differential equation (PDE) frameworks have combined evaporation and fluorescein dynamics to estimate otherwise unmeasurable parameters \citep{zhongDynamicsFluorescentImaging2019,lukeParameterEstimation2020,lukeParameterEstimationMixedMechanism2021}. More recently, spatially independent ODE models have been calibrated to small TBU spots and streaks and coupled with neural-network based data extraction to expand datasets \citep{lukeFittingSimplifiedModels2021,driscollFittingODEModels2023}. Driscoll et al. combined reduced ODE models with a convolutional-neural-network-based image analysis pipeline that automatically identified multiple tear-breakup regions in fluorescein video sequences and extracted their intensity time series. This automation substantially increased the number of analyzable breakup events, enabling model fitting and parameter estimation for hundreds of TBU instances across multiple subjects.

Additional theoretical and experimental studies have investigated mechanisms underlying thin-film rupture and instability. \citet{ji2016finitetimerupture} analyzed a lubrication-theory thin-film model with modified evaporative loss and used self-similar asymptotic analysis together with high-resolution adaptive finite-difference simulations to characterize evaporation-driven finite-time rupture singularities. \citet{shi2022instabilitybubble} combined lubrication-theory modeling and linear stability analysis to show how interactions between Marangoni and capillary flows can induce symmetry breaking and instability in surfactant-laden thin films.

While the results of \citet{driscollFittingODEModels2023} are an important first step toward a complete understanding of TF dynamics, they are limited on the one hand by the relative crudeness of the mathematical model and on the other hand by the need to computationally solve an inverse problem for each prospective TBU location. While fully 2D PDE model simulation has been established \citep{chenEvaporationdrivenTearFilm2024} and simple inverse problems for the 2D model are underway for isolated TBU cases \citep{ChenDriscoll2D}, solving the 2D inverse problem at large scale is currently computationally expensive. 

Recent advances in operator learning have produced a variety of architectures for approximating mappings between function spaces. \citep{Tancik2020} showed that standard multilayer perceptrons exhibit spectral bias toward low-frequency functions and proposed Fourier feature mappings to improve the learning of high-frequency components, significantly enhancing the performance of coordinate-based neural networks on low-dimensional regression tasks. \citet{LiFourierNeural2021a} introduced the physics-informed neural operator (PINO), which incorporates PDE constraints into neural operator training to improve accuracy and generalization, particularly when training data are limited.  \citet{WangEigenvectorBias2021} analyzed the training dynamics of physics-informed neural networks (PINNs) for multiscale PDEs and proposed Fourier feature embeddings to improve their representation of high-frequency solutions. \citet{ChenOperatorLearning2023} combined Fourier neural operators with neural ODEs to develop continuous spatial-temporal operator learning that is invariant to both spatial and temporal discretizations and improves long-term prediction stability. \citet{JMLR:v24:23-0478} established a rigorous theoretical foundation for PCA-Net by proving universal approximation results and deriving complexity bounds for PCA-based operator learning in infinite-dimensional spaces.  

Operator learning has also been successfully applied to a wide range of scientific computing problems. \citet{Lu2021DeepONet} introduced DeepONet and demonstrated its effectiveness for learning solution operators of nonlinear ODEs and PDEs. \citet{Kovachki2021NeuralOperator} showed that neural operators accurately approximate solution operators for benchmark problems including Burgers' equation, Darcy flow, and the Navier–Stokes equations with substantial computational speedups over traditional numerical solvers. \citet{Lu2022MultifidelityDeepONet} further developed Multifidelity DeepONet, which combines low- and high-fidelity data to efficiently learn nonlinear operators and demonstrated its effectiveness on several forward and inverse PDE problems.

Motivated by these recent advances, we investigate operator learning techniques to ODE and 1D PDE models of TBU dynamics. In operator learning, one trains a deep neural network or other learning architecture to learn the mapping implied by the solution of a differential equation. This could mean the mapping from parameters to solutions, for instance, but our goal is to learn the mapping from observations, in the form of synthetic FL imaging data, to outcomes, in the form of TF thickness and osmolarity functions, without the need to simulate the physics or solve an inverse problem. If successful, this process provides a fast computational method that could be used in real time to predict characteristics of tear film flows from FL imaging. At the modeling level, we can compare the results of predictors trained on ODE and PDE models to clarify how often and under what circumstances they make substantially different predictions.

In \autoref{sec:mathmodels} we describe the ODE and PDE models used to simulate TF dynamics during TBU. In \autoref{sec:operator-learning} we describe three operator learning methods that are used and tested. In \autoref{sec:results} we present the results of numerical experiments with the two simulation models and operator learning methods, and in \autoref{sec:discussion} we make concluding remarks and observations.

\section{Mathematical models}
\label{sec:mathmodels}

\subsection{One-dimensional PDE}

\begin{figure}
    \centering
    \includegraphics[width=0.8\textwidth]{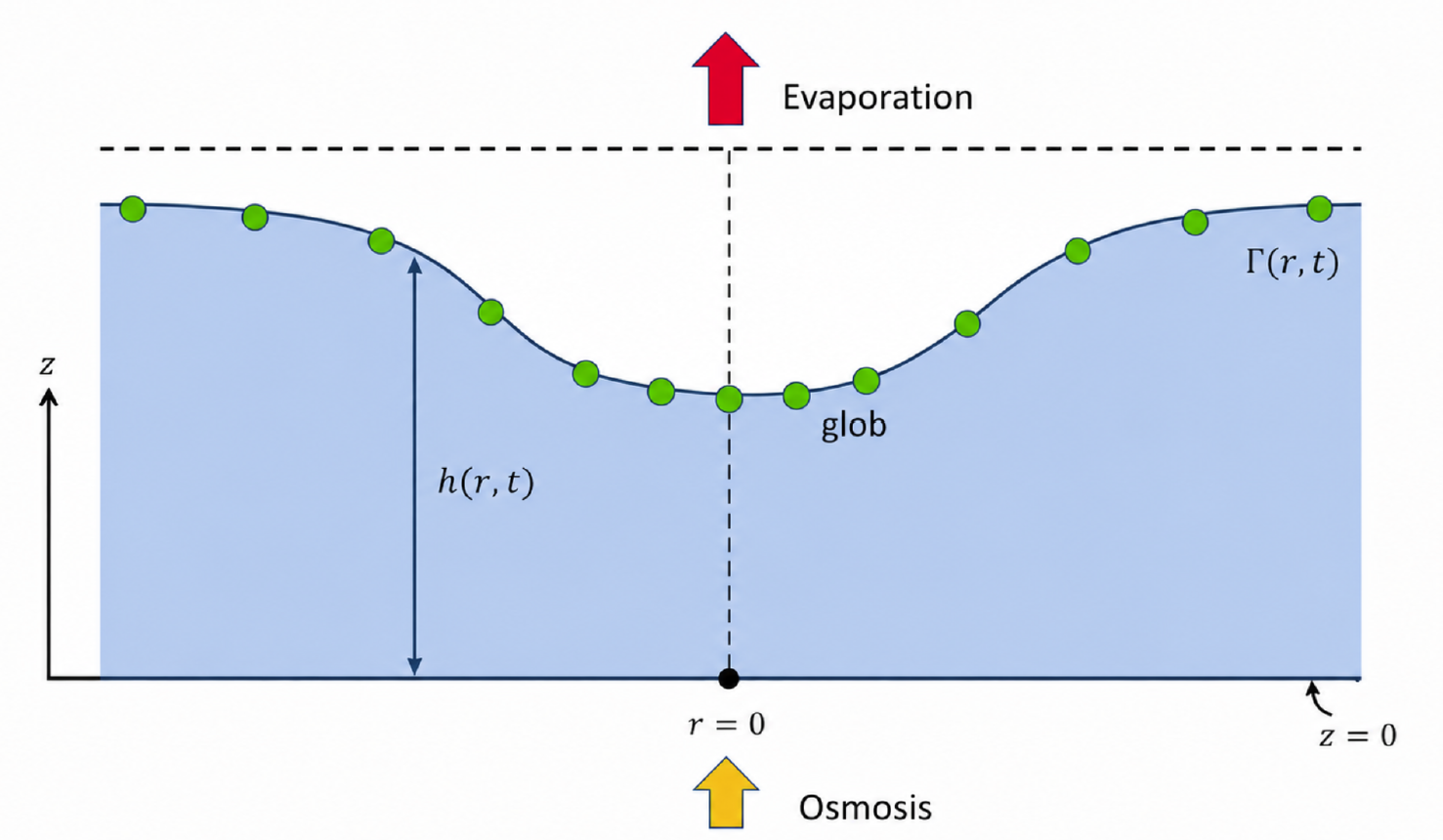}
    \caption{Mechanisms in the tear film at a later time $t>0$. Evaporation occurs at the top, osmosis at the bottom; The green globs represent localized excess lipids at the air–tear interface. $\Gamma$ represents the surfactant concentration, and $h$ represents the tear film thickness. Here, $r=0$ denotes the axis of symmetry, and $z=0$ corresponds to the corneal surface.}
    \label{fig:pde_mec}
\end{figure}

The PDE system is derived from \citet{zhongDynamicsFluorescentImaging2019} and describes the situation depicted in \autoref{fig:pde_mec} of a radially symmetric glob of lipid on top of the aqueous layer at a later time $t>0$. The dashed line in the picture represents the tear film at $t=0$. We assume there is nonzero constant evaporation under the glob and zero evaporation outside the glob. The Marangoni effect, driven by tangential flow, arises in regions where thicker globs of the lipid layer are present. These areas, with higher surfactant concentration, reduce the aqueous/air surface tension and trigger fluid motion along the surface, leading to rapid film breakup. We expect that diffusion, advection, evaporation, osmosis, and Marangoni flow together influence the TF thickness which would result in TBU.  

The independent and dependent variables are shown in \autoref{tab:variables}. Quantities named without primes are nondimensionalized according to 
\begin{gather*}
r'=\ell r , \quad z'=h_0'z, \quad t'=t_s't, \quad h'=h_0'h,\quad u_r'=Uu_r,  \\
\quad p'=\frac{\mu U}{\ell \epsilon^2}p, \quad \Gamma'=\Gamma_0\Gamma, \quad J'=\epsilon\rho UJ, \quad c'=c_0 c,\quad f'=f_{\text{cr}}f, \quad f_0'=f_{\text{cr}}f_0
\end{gather*}
where primes denote dimensional quantities that are eventually replaced by nondimensionalized versions. The relevant physical parameters are given in \autoref{tab:physical}, and the nondimensional parameters are shown in \autoref{tab:nondim_pde}. 

\begin{table}[tbp]
    \centering
    \begin{tabular}{cl}
        Variable & Meaning \\ \hline
        $r$ & Spot radius \\
        $z$ & Depth dimension \\
        $t$ & Time \\
        $\Gamma(r,t)$ & Lipid surfactant concentration \\
        $u_r(r,h,t)$ & Horizontal surface fluid velocity \\
        $h(r,t)$ & TF thickness \\ 
        $p(r,t)$ & Pressure \\
        $J(r)$ & Evaporation rate \\ 
        $c(r,t)$ & Osmolarity \\ 
        $f(r,t)$ & Fluorescein concentration
    \end{tabular}
    \caption{Variables in the 1D glob model. }
    \label{tab:variables}
\end{table}

\begin{table}[tbp]
    \centering
    \begin{tabular}{c l c l} 
      Parameter & Description & Value & Reference \\ 
     \hline
     $\mu$ & Viscosity & \qty{1.3e-3}{\Pa\second} & \citet{tiffanyViscosityHumanTears1991} \\
     $\sigma_{0}$ & Surface tension & \qty{0.045}{\N\per\meter} & \citet{nagyovaComponentsResponsibleSurface1999} \\
     $\rho$& Density & \qty{1000}{\kg\per \meter\tothe{3}} & Water \\
     $A^*$ & Hamaker constant & $6\pi \cdot 3.5 \times 10^{-19}$ J & \citet{braunTearFilm2018} \\
     $h_0$ & Initial TF thickness & \qtyrange{2}{8}{\micro\meter} & Calculated \\  
     $(\Delta\sigma)_0$ & Change in surface tension & \qtyrange{2}{60}{\micro\N\per\meter} & Calculated \\ 
     $t_s$ & Observation timescale & \qtyrange{10}{60}{\second} & Calculated \\ 
     $v$ & Thinning rate & \qtyrange{0}{40}{\micro\meter\per\minute} & Calculated \\ 
    $f_0$ & Initial FL concentration & \qtyrange{0.08}{0.2}{\%} & Calculated \\ 
    $R_I$ & Glob radius & \qtyrange{0.02}{0.2}{mm} & Calculated \\ 
     $V_{w}$ & Molar volume of water & \qty{1.8e-5}{\meter\cubed \per \mol} & Water \\ 
     $D_{f}$ & Diffusivity of fluorescein & \qty{0.39e-9}{\meter\squared \per \second} &  \citet{casaliniDiffusionAggregationSodium2011} \\ 
     $D_{o}$& Diffusivity of salt & \qty{1.6e-9}{\meter\squared \per \second} & \citet{riquelmeInterferometricMeasurementDiffusion2007} \\
     $c_{0}$& Isotonic osmolarity & \qty{300}{\Osm \per \meter\cubed} &  \citet{lempTearOsmolarityDiagnosis2011}\\
     $P_{0}$& Permeability of cornea & \qty{12.1}{\um\per\second} &  \citet{braunDynamicsFunctionTear2015}\\
     $\ell$ & Characteristic length, $\left({t_s\sigma_0h_0^3}/{\mu}\right)^{\frac{1}{4}}$ & \qtyrange{0.138}{0.412}{mm} & Calculated \\ 
     $U$ & Characteristic velocity, $\left({\sigma_0h_0^{3}}/{\mu t_s^3}\right)^{\frac{1}{4}}$ & \qtyrange{0.056}{0.099}{\milli\meter\per\second} & Calculated \\ 
     $\epsilon_{f}$& Napierian extinction coefficient & \qty{1.75e7}{L / m mol} & \citet{motaSpectrophotometricAnalysisSodium1991}\\
     $f_{cr}$& Critical FL concentration  &  $0.2\%$ & 
     \citet{webberContinuousFluorophotometricMethod1986}
    \end{tabular}
    \caption{Physical parameters (dimensional) used in the governing equations.}
    \label{tab:physical}
\end{table}

\begin{table}[tbp]
    \centering
    \begin{tabular}{c l c}
        Parameter & Description & Expression \\
        \hline
        $\epsilon$ & Aspect ratio & $h_0' / \ell$ \\
        $M$ & Marangoni number & $\epsilon (\Delta\sigma)_0 \left[ t_s'^3 / (\sigma_0 \mu^3 h_0'^3) \right]^{1/4}$ \\
        $A$ & Nondimensional Hamaker constant & $A^* / \left[ \epsilon (\Delta\sigma)_0 h_0' \ell \right]$ \\
        $P_c$ & Permeability of cornea & $(P_0 V_w c_0) / (\epsilon U)$ \\
        $\mathrm{Pe}_f$ & P\'{e}clet number for FL diffusion & $U \ell / D_f$ \\
        $\mathrm{Pe}_c$ & P\'{e}clet number for salt ion diffusion & $U \ell / D_o$ \\
        $\mathrm{Pe}_s$ & P\'{e}clet number for surface diffusion & $\epsilon (\Delta\sigma)_0 \ell / (\mu D_s)$ \\
        $\phi$ & Nondimensional Napierian extinction coefficient & $\epsilon f_{cr} h_0'$ \\
    \end{tabular}
    \caption{Nondimensional parameters used in the governing equations.}
    \label{tab:nondim_pde}  
\end{table}

The resulting nondimensional PDE system for the dependent variables $h$, $p$, $c$, $\Gamma$ and $f$ is
\begin{subequations}
\label{eq:pde1d}
    \begin{align}
    \partial_t h &= -J + P_c(c - 1) - \frac{1}{r} \partial_r (r h \overline{u}), \label{eq:pdeh}\\
    p &= -\frac{1}{r} \partial_r (r \partial_r h) - A h^{-3}, \label{eq:pdep}\\
    \partial_t \Gamma &= \left[\text{Pe}_s^{-1} \left( \frac{1}{r} \partial_r (r \partial_r \Gamma)\right) - \frac{1}{r} \partial_r (r u_r \Gamma) \right]B, \label{eq:pdegamma}\\
    h(\partial_t c + \overline{u} \partial_r c) &= \text{Pe}_c^{-1} \left( \frac{1}{r} \partial_r (r h \partial_r c) \right) + Jc - P_c(c -1)c\\
    h(\partial_t f + \overline{u} \partial_r f) &= \text{Pe}_f^{-1} \left( \frac{1}{r} \partial_r (r h \partial_r f) \right) + Jf - P_c(c -1)f, \label{eq:pdef}
    \end{align}
\end{subequations}
where $B$ is a smooth approximation to a transition step function,
\[
B(r,R_I) = \frac{1}{2} + \frac{1}{2} \tanh \left( \frac{r - R_I}{0.1} \right),
\]
and the horizontal surface TF fluid velocity \( u_r(r,h,t) \) and average horizontal TF fluid velocity \( \overline{u} \) throughout the film are given by
\begin{align}
u_r(r, h, t) &= -\frac{\frac{1}{2}h^2 (\partial_r p) B + M \partial_r \Gamma h B}{B+(1-B)h}, \\
\overline{u} &= -\frac{\frac{1}{3}h^2 (\partial_r p)[B+\frac{1}{4}h(1-B)] + \frac{1}{2}M \partial_r \Gamma h B}{B+(1-B)h}.
\end{align}

The governing equations are obtained from the incompressible Navier--Stokes equations under the lubrication approximation, using the small aspect ratio of the tear film to reduce the momentum equations to a balance between viscous stresses, capillary pressure gradients, and Marangoni stresses. The pressure is independent of depth at leading order and is determined by the curvature of the free surface,
\[
p=-\frac{1}{r}\partial_r(r\partial_r h),
\]
Conservation of mass over the film thickness gives the evolution equation
\[
\partial_t h=-J+P_c(c-1)-\frac{1}{r}\partial_r(rh\overline{u}),
\]
in which $-J$ represents loss of water by evaporation, $P_c(c-1)$ represents osmotic flow from the corneal surface into the tear film when the osmolarity is above isotonic, and the divergence term describes redistribution of tear fluid by the depth-averaged radial velocity $\overline{u}$. Thus, localized evaporation and divergent Marangoni-driven flow thin the film, while osmosis and capillary-driven flow may partially oppose thinning.
In order to simulate the effect of the lipid glob on evaporation, we define the spatial dependence of the evaporative flux by
\[
J(r) = (1 - B(r)) v.
\]

The depth-averaged velocity $\overline{u}$ and the surface velocity $u_r$ are obtained by solving the leading-order radial momentum equation with no slip at the corneal surface and a mixed tangential stress condition at the tear--air interface. The pressure-gradient terms in these velocities represent capillary-pressure-driven flow, while the terms proportional to $M\partial_r\Gamma$ represent Marangoni flow generated by gradients in surfactant concentration. The surfactant equation
\[
\partial_t\Gamma =
\left[
\mathrm{Pe}_s^{-1}\frac{1}{r}\partial_r(r\partial_r\Gamma)
-\frac{1}{r}\partial_r(ru_r\Gamma)
\right]B
\]
describes surface diffusion and advection of lipid surfactant outside the glob region. The factor $B(r)$ smoothly switches the surfactant dynamics and surface mobility between the immobile glob-covered region and the mobile exposed region. Finally, the osmolarity and fluorescein equations are depth-averaged advection--diffusion equations. In both equations, the left-hand side represents local change and advection by $\overline{u}$, the first term on the right-hand side represents diffusion, the $Jc$ or $Jf$ term represents concentration increase due to evaporative water loss, and the $-P_c(c-1)c$ or $-P_c(c-1)f$ term represents dilution by osmotic inflow.

Since we are not interested in boundary effects such as eyelids, we assume periodic conditions on all the dependent variables and that all the dependent variables are initially uniform:
\begin{equation}
\label{eq:ic_pde}
    h(r,0) = c(r,0) = 1,\quad f(r,0) = f_0,
\end{equation}
where $f_0$ is the FL concentration normalized to the critical concentration $f_{\text{cr}}$.
Once $h$ and $f$ are known, we compute the FL intensity via~\citep{webberContinuousFluorophotometricMethod1986,braunModelTearFilm2014}
\begin{equation}
\label{eq:FLI}
I=I_0\frac{1-\exp(-\phi fh)}{1+f^2}, 
\end{equation}
where $I_0$ is a normalization coefficient so that $I(0)=1$ and $\phi$ is the nondimensional Napierian extinction coefficient.

The values of $h_0'$, $t_s'$, and $f_0'$ are determined for each FL trial by the methods described in \citet{wuEffectsIncreasingOcular2015}, and we consider them as being provided to the model. The glob radius $R_I'$, surface tension change $(\Delta\sigma)_0'$, and nominal thinning rate $v'$ are additional parameters to the model. Hence, there is a 6-dimensional parameter space affecting the output of the model. 

For discretization, we use the method of lines with a Fourier spectral collocation method~\citep{trefethenSpectralMethodsMatlab2000} in space on a uniform periodic grid and exploit symmetry by solving over $r \in(0,\pi]$. The resulting discretization of spatial terms in \eqref{eq:pdeh}--\eqref{eq:pdef} creates a differential--algebraic system (DAE) that is solved in Julia using the Rodas5P solver, a 5th order A-stable stiffly stable Rosenbrock method in the \textbf{DifferentialEquations} package~\citep{rackauckas2017differentialequations}. Computations were performed using Julia version $1.10.5$ on a Windows laptop with a 12th Gen Intel Core i7-12700H 2.30GHz processor and 32GB RAM.  

\subsection{ODE}

To simplify the model \eqref{eq:pde1d}, we follow the procedure in~\citet{driscollFittingODEModels2023} to remove the spatial dependence and model only at the glob center. In order to replace the tangential flow caused by surface tension gradients that are no longer possible, we impose a time-dependent shear stress term
\begin{equation}
    \label{eq:shear}
    g(t) = b_1' e^{-b_2't},
\end{equation}
for dimensional parameters $b_1'$, the maximum shear, and $b_2'\ge 0$, the decay constant. When $b_1'>0$, the film thins as though there is outward flow from the center of the glob, while $b_1'<0$ causes thickening from an inward flow. Because there can no longer be diffusion, the only difference between osmolarity $c(t)$ and fluorescein concentration $f(t)$ is a constant of proportionality, so $f$ is dropped as a dependent variable.
In addition, we assume a constant evaporation rate $J_e'$.
Then we nondimensionalize the parameters and variables according to the following relations:
\[
\begin{aligned}
& J_e = \frac{J_e'   t_s'}{h_0'}, \quad f_0 = \frac{f_0'}{f_{\text{cr}}}, \quad b_1 = b_1'   t_s', \quad b_2 = b_2'   t_s'
\end{aligned}
\]
where $f_{\text{cr}} = 0.2$ \% is the critical fluorescence concentration.

The resulting system, in nondimensional form, is 
\begin{subequations}
    \label{eq:ode}
    \begin{align}
    \frac{dh}{dt} &= P_c(c-1) - g(t)h - J, \label{eq:odeh}\\
    h\frac{dc}{dt} &= - g(t)hc - \frac{dh}{dt}c, \label{eq:odec}\\
    h(0) &= 1,
    \label{eq:ich}\\
    c(0) &= 1.
    \label{eq:icc}
    \end{align}
\end{subequations}
This model was labeled as Model D in~\citet{driscollFittingODEModels2023}. Once $h$ and $c$ are computed, $f(t)$ is found via $f(t) = f_0\, c(t)$. The FL intensity $I$ is again given by~\eqref{eq:FLI}.

The dimensional parameters for~\eqref{eq:ode} again include the initial thickness $h_0'$, the trial time $t_s'$, the initial fluorescein concentration, $f_0'$, and the thinning rate $J_e'$, as in the PDE model. However, the glob radius and surface tension parameters are replaced by the constants $b_1'$ and $b_2'$ in~\eqref{eq:shear}.

\section{Operator learning learners}
\label{sec:operator-learning}

Each of the mathematical models in \autoref{sec:mathmodels} effectively maps a parameter space to physically relevant outputs for nondimensional (i.e., relative) thickness $h$, osmolarity $c$, and the observed FL intensity $I$. In practice, our goal is to map an observed time series for $I(t)$ to the thickness and osmolarity functions. Using a mathematical model, this can be treated as an inverse problem, using an observation to find the best-fit parameters of the model, whose solution provides $h$ and $c$.

Instead, we can use operator learning~\citep{ChenOperatorLearning2023} to train a machine learning method to produce the $I\mapsto(h,c)$ mapping directly, without solving or even referring to a mathematical model. This approach offers the possibility of using model solutions solely to produce training data. After the training is complete, the learner can make predictions very quickly. In addition, by training learners on different mathematical models, we can assess to what extent, and for which observations, one model produces results similar to the other.

There are many forms of operator learning in the literature, using different ways to represent the input and output spaces as well as the connection between them. We report results for three approaches: Fourier feature networks (FFN), dense neural networks between spaces that are compressed via principal components analysis (Dense-PCA), and dense networks with inputs that enhance PCA inputs with externally known parameters (Dense-PCAX). 

\subsection{Generation of training/testing data}

The ODE and PDE models of \autoref{sec:mathmodels} were used to generate two synthetic datasets for training and testing the ML learners. Each begins by sampling parameters using six-dimensional Halton sequences within the ranges given in \autoref{tab:paramranges}. Because it is not straightforward to \textit{a priori} characterize which parameters produce physically and physiologically reasonable results, we reject simulations that have any of the following characteristics:
\begin{enumerate}
    \item \textbf{Upper bound on thickness.} Instances where $h(t) > 1.1$ at any time were discarded, as tear film thickness is generally expected to decrease at the start of simulation time.
    \item \textbf{Lower bound on thickness.} Solutions where $h(t) < 0.2$ were excluded, since different dynamics will hold when the thickness of the film is comparable to the glycocalyx thickness.
    \item \textbf{Nonphysical growth.} Instances with $h(1) > 1.5  \min h(t)$ were removed to rule out nonphysical regrowth after an initial thinning phase.
    Here $t=1$ is the nondimensional final time for all cases.
    \item \textbf{Increasing intensity.} We collect data only with decreasing fluorescence intensity over time, so instances where $I(t)$ increases were discarded.
\end{enumerate}
The surviving cases constitute a parameter set $\Omega \subset \mathbf{R}^6$. 

\begin{table}[tbp]
    \centering
    \begin{tabular}{cc|cc}
       \multicolumn{2}{c|}{ODE}  &  \multicolumn{2}{c}{PDE}   \\
        Parameter & Range & Parameter & Range \\ \hline
        $h_0'$ & $[\si{2}{\micro\meter}, \si{8}{\micro\meter}]$ &
          $h_0'$ & $[\si{2}{\micro\meter}, \si{8}{\micro\meter}]$ \\
        $f_0'$ & $[\si{0.08}{\%}, \si{0.2}{\%}]$ &
          $f_0'$ & $[\si{0.08}{\%}, \si{0.2}{\%}]$  \\
        $t_s'$ & $[\si{10}{\second}, \si{60}{\second}]$ &
          $t_s'$ & $[\si{10}{\second}, \si{60}{\second}]$ \\
        $v'$ & $[\si{0}{\micro\meter\per\minute}, \si{40}{\micro\meter\per\minute}]$ &
         $J_e'$ & $[\si{0}{\micro\meter\per\minute}, \si{40}{\micro\meter\per\minute}]$ \\
        $R_I'$ & $[\si{0.02}{\milli\meter},\si{0.2}{\milli\meter}]$ &
          $b_1$ & $[\si{-0.2}{\second^{-1}}, \si{2.0}{\second^{-1}}]$ \\
        $(\Delta\sigma)_0'$ & $[\si{2}{\micro\N\per\meter},\si{60}{\micro\N\per\meter}]$ &
          $b_2$ & $[\si{0}{\second^{-1}}, \si{2.0}{\second^{-1}}]$ 
    \end{tabular}
    \caption{Ranges for dimensional parameters used to create synthetic datasets for ML training. See \autoref{sec:mathmodels} for their roles in the mathematical models.}
    \label{tab:paramranges}
\end{table}

Each element of $\Omega$ is associated with resulting solutions for $c$, $h$, and $I$. For the ODE, each variable is a function of time only and was sampled at $N=601$ equally spaced times; for the PDE, each was sampled at $N$ equally spaced times and $r=0$, which is the center of the glob. These discretizations produced the realization sets $\mathcal{C}$, $\mathcal{H}$, and $\mathcal{I}$, all subsets of $\mathbf{R}^N$. The elements of $\mathcal{I}$ were used as input feature vectors mapping to corresponding output feature vectors from either $\mathcal{C}$ or $\mathcal{H}$. (It would also be possible to concatenate corresponding elements of $\mathcal{C}$ and $\mathcal{H}$ and train a single learner to produce both outputs simultaneously, but that possibility was not explored for this work.) For the PDE model, $|\mathcal{I}| = 28698$, while for the ODE model, $|\mathcal{I}| = 64522$.

In order to improve the robustness of predictions for noisy experimental observations, we augmented the input set with copies that were randomly perturbed by multiplicative noise. To make perturbations similar to those observed in the experimental dataset, the noise for a signal vector $I$ is a Gaussian-smoothed (\(\sigma = 9\)) white noise of standard deviation $0.25 |I_1 - I_N|$. Two perturbed copies of each input vector (and clean copies of the corresponding outputs) were added to $\mathcal{I}$, which was then split into training and testing subsets \train{I} (75\%) and \test{I} (25\%), respectively, along with the corresponding outputs.  

\subsection{Fourier feature network (FFN)}

A Fourier feature network (FFN) uses an input layer of randomized linear combinations of Fourier modes of the input features to reduce dimensionality while retaining the ability to learn functions with fine-scale variations~\citep{Wang21}. Our FFN comprises an input layer, three hidden layers with ReLU activation, and an output layer:
\[
\hat{y} = f(\gamma(I)) = W_4 \, \sigma(W_3 \, \sigma(W_2 \, \sigma(W_1 \gamma(I) + v_1) + v_2) + v_3) + v_4.
\]
Here, $\gamma(I)$ denotes the Fourier feature mapping, defined as
\[
\gamma(I) = [\sin(2\pi \sigma_j I), \cos(2\pi \sigma_j I)]_{j=1}^{3},
\]
where $I$ is the $N$-dimensional input vector and $\sigma_j \in \{0.5, 1, 2\}$ represents frequency scales selected to capture varying resolution details. This makes $\gamma(I)$ an array of length $6N$. The hidden layers are dense networks with weight matrices $W_i$ and bias vectors $v_i$ for $i=1,2,3$, and $W_4$ and $v_4$ characterize the dense output layer to produce a time series prediction of length $N$ for either $c$ or $h$. 

The first weight matrix $W_1$ has shape $300\times 6N$ and the bias vector $v_1$ has length $300$. The next two hidden layers have weight matrices of sizes $200\times300$ and $200\times200$, with corresponding bias vectors. After each of the first three layers, the activation function $\sigma=\mathrm{ReLU}$ is applied.
The last weight matrix $W_4$ of the output layer has size $N\times200$ to produce a time series of either $c$ or $h$.

Training was performed by solving the optimization problem
\begin{equation}
\min_{\theta}\;
\mathcal{L}(\theta)
=
\frac{1}{N}\sum_{t=1}^{N}\left(y_t-\hat{y}_t(\theta)\right)^2
+
\frac{\alpha}{N-1}\sum_{t=1}^{N-1}
\left(\hat{y}_{t+1}(\theta)-\hat{y}_t(\theta)\right)^2,
\label{eq:training_loss}
\end{equation}
where $\theta$ denotes all trainable network parameters (weights and biases), and
$\hat{y}(\theta)$ is the network prediction. The second term suppresses
high-frequency oscillations, since $\hat{y}$ represents a time series from a
dynamical system with smooth solutions. We chose $\alpha=1.6$ as the regularization parameter.

\subsection{Dense neural network with PCA (Dense-PCA)}

Another way to encode the input and output time series is to apply PCA to them, connecting the coefficient spaces by a fully-connected neural network~\citep{Yadav19}. Reducing dimensionality reduces the number of training parameters and may help focus attention on the most salient features of the signals.

 The fully connected neural networks used three hidden layers of 300 neurons with ReLU activation:
\begin{align}    
\label{eq:denselayers}
    x &= P(I-\mu)\\
   \hat{y} &= \nu + Q [W_4 \, \sigma(W_3 \, \sigma(W_2 \, \sigma(W_1 x + v_1) + v_2) + v_3) + v_4],
\end{align}
where $(P,\mu)$ and $(Q,\nu)$ denote PCA compression and reconstruction, respectively, and the $W_i$ and $v_i$ are the weights and biases of the dense hidden layers, determined through minimizing the MSE loss between the training samples and their corresponding predictions.

\subsection{Dense network on PCA augmented with external parameters (Dense-PCAX)}

In \citet{driscollFittingODEModels2023}, the values of the initial thickness $h_0'$, initial fluorescein concentration $f_0'$, and observation time length $t_s'$ are all determined by processes external to the inverse problem for the ODE, which optimized for the remaining parameters describing the evaporation rate and flow. 

Our ODE and PDE models produce training data by varying over {both the shared dimensional parameters and the model-specific parameters} without making a distinction. In order to determine whether forecasting ability improves when given access to exactly the same data as the ODE inverse problem in~\citet{driscollFittingODEModels2023}, we ran a variant of Dense-PCA in which the PCA component vector $x$ in~\eqref{eq:denselayers} was concatenated with the values of $h_0'$, $f_0'$, and $t_s'$, which are included as inputs to the learner. 

\section{Results}
\label{sec:results}

\begin{figure}
    \centering
    \includegraphics[width=0.9\textwidth]{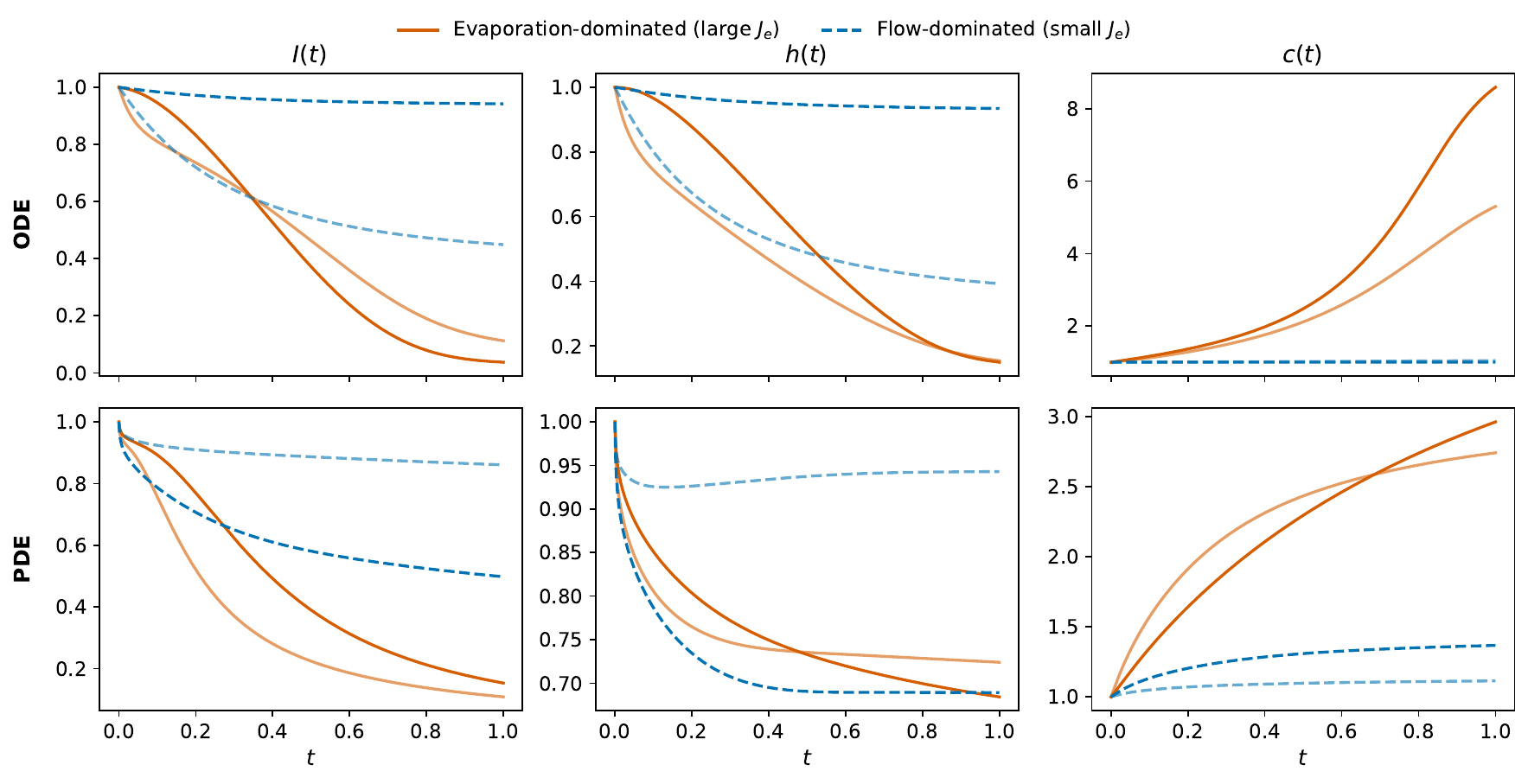}
    \caption{Selected samples from the synthetic datasets generated by the ODE and PDE models. Each row shows a different sample, with the left column showing the fluorescein intensity $I(t)$, the middle column showing the thickness $h(t)$, and the right column showing the osmolarity $c(t)$. The top three rows are from the ODE model, while the bottom three rows are from the PDE model.}
    \label{fig:synthetic-samples}
\end{figure}

As described in \autoref{sec:operator-learning}, the two mathematical models (\autoref{eq:pde1d} and \autoref{eq:ode}) each generated a synthetic dataset for three ML learners. Selected samples from these datasets are shown in \autoref{fig:synthetic-samples}.

For the input training set \train{I}, 4 PCA components were sufficient to achieve reconstruction mean squared error (MSE) of approximately $1.4 \times 10^{-4}$. For the output training sets, the number of components needed for comparable reconstruction varied by mathematical model and output variable. The sizes of the layers in the ML learners are summarized in \autoref{tab:model-params}. Because these models are all small by modern standards, the learners were trained for a relatively long time (1000 epochs using the Adam optimizer with a learning rate of 0.001) in order to reach optimal weights.

\begin{table}[tbp]
    \centering
    \begin{tabular}{l|c|c|c}
        \textbf{Learner} & \textbf{Input} & \textbf{Output} & \textbf{Hidden layers} \\
        %& \textbf{Trainable parameters}\\
        \hline   
        FFN &
        601 & 601 & 300, 200, 200 \\
        % & 401801\\
        Dense-PCA &
        4 & 4 (ODE) / 13 (PDE, $c$) / 8 (PDE, $h$) & 300, 300, 300 \\ 
        %& 183304  \\
        Dense-PCAX &
        7 & 4 (ODE) / 13 (PDE, $c$) / 8 (PDE, $h$) & 300, 300, 300 \\
        % & 184204  \\
    \end{tabular}
    \caption{Neural network layer sizes used for our results.}
    \label{tab:model-params}
\end{table}

\begin{figure}
    \centering
    \includegraphics[width=1.0\textwidth]{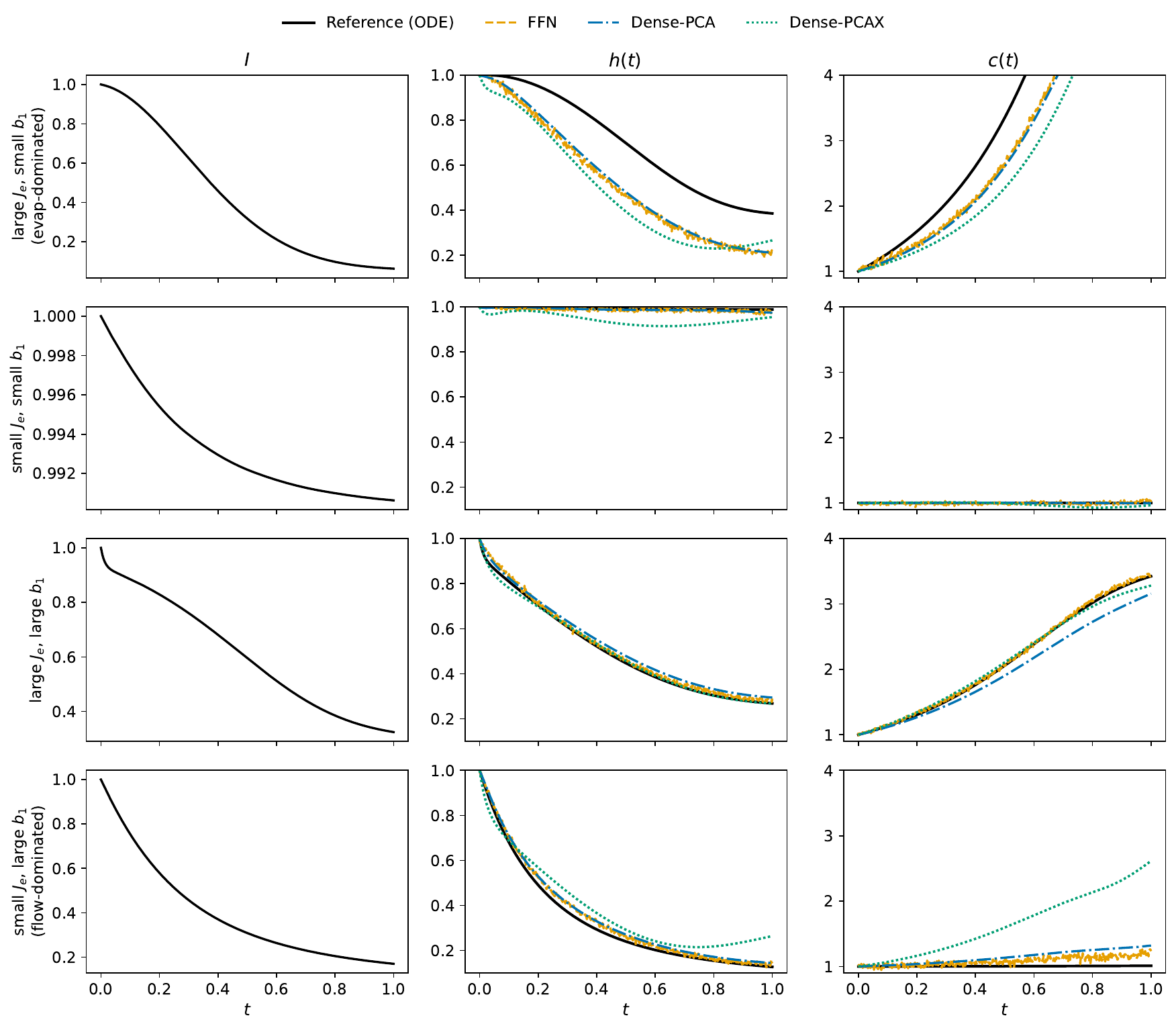}
    \caption{Example predictions on synthetic data by the three operator learners trained on ODE data.}
    \label{fig:example-predictions-synthetic-ode}
\end{figure}

\begin{figure}
    \centering
    \includegraphics[width=1.0\textwidth]{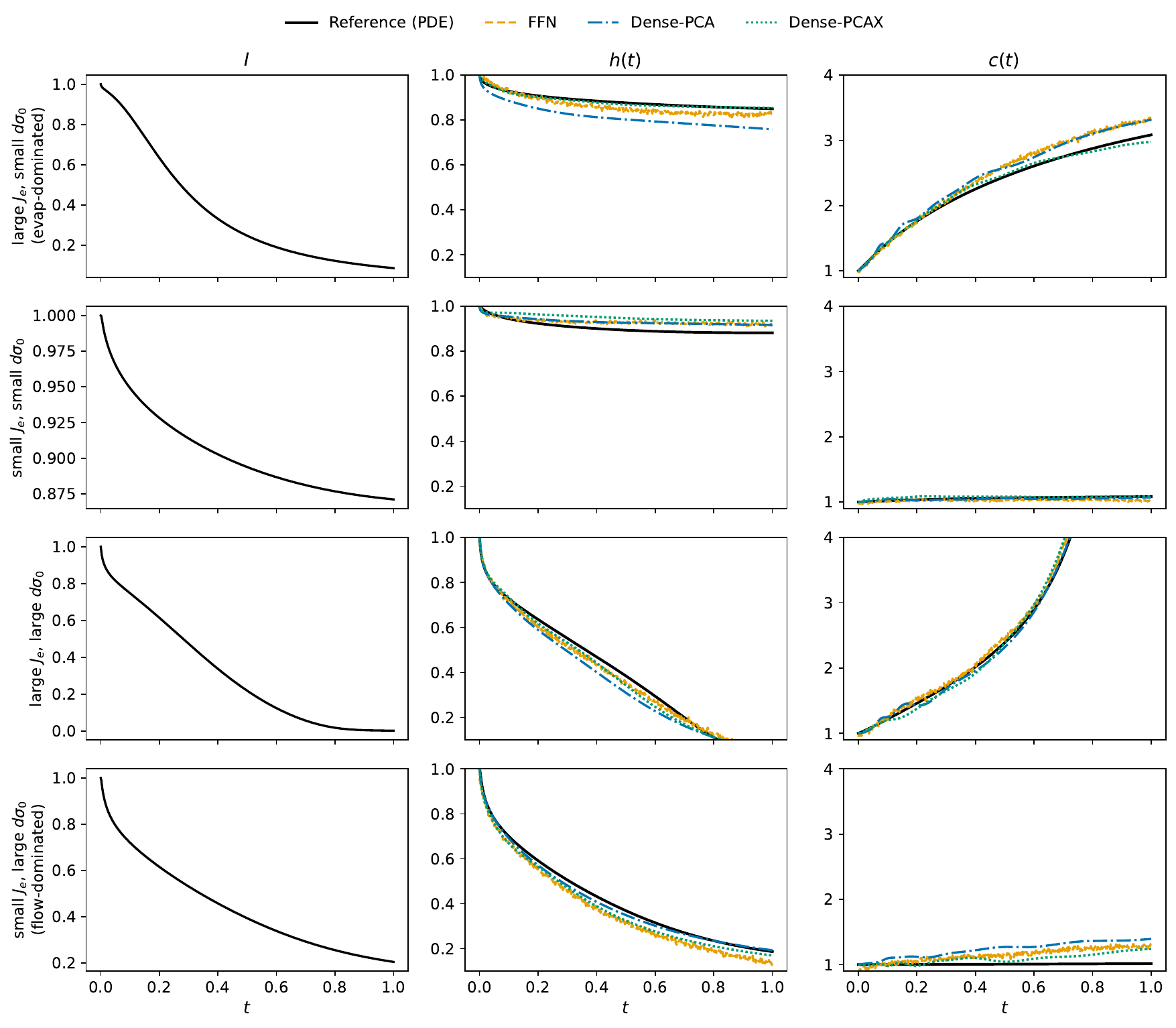}
    \caption{Example predictions on synthetic data by the three operator learners trained on PDE data.}
    \label{fig:example-predictions-synthetic-pde}
\end{figure}

\autoref{fig:example-predictions-synthetic-ode} shows example predictions of the three learners trained on ODE data, while \autoref{fig:example-predictions-synthetic-pde} shows them for learners trained on PDE data. These figures show only unperturbed intensity inputs; noisy corruptions of them are also part of the testing set. Across all the cases, the FFN predictions show a great deal of high-frequency oscillation, which is not characteristic of the true solutions, but these predictions can accurately follow the true curves. The PCA-based learners always produce smooth predictions. 

To summarize our results over the full training sets, we report distributions of relative root-mean-square error (rRMSE), defined for a prediction $\hat{y}$ of time series vector $y$ by
\begin{align*}
\mathrm{RMSE}_y &=\left[ \tfrac{1}{N} \sum_{k=1}^{N} \left(y_k- \hat{y}_k \right)^2\right]^{1/2},\\
\mathrm{rRMSE}_y  &= \frac{\mathrm{RMSE}_y}{\left(\tfrac{1}{N}\sum_{k=1}^{N} y_k^2\right)^{1/2}}.
\end{align*}

\subsection{Synthetic data}

\begin{figure}
  \centering
  \includegraphics[width=1.0\textwidth]{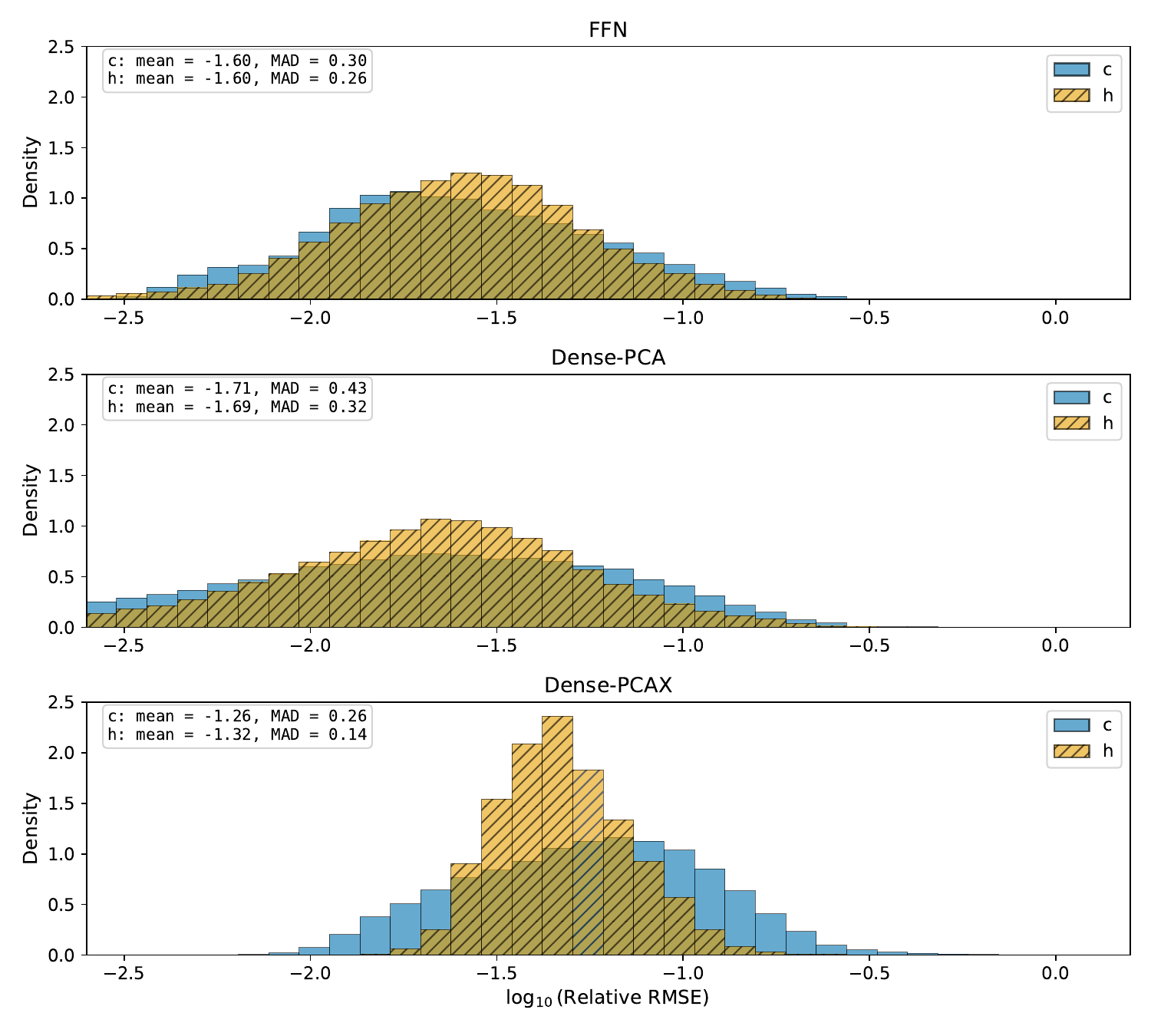}
  \caption{Relative RMSE for synthetic testing of ML learners trained by ODE data. The mean and mean absolute deviation (MAD) is shown for each distribution. Top to bottom: FFN, Dense-PCA, Dense-PCAX.}
  \label{fig:ode-ffn-hists}
\end{figure}

\autoref{fig:ode-ffn-hists} shows the distributions of rRMSE of both $h$ and $c$ for all three ODE-trained learners on the synthetic test set. Thickness $h$ was predicted to at least one digit of accuracy in nearly all cases, as was $c$ for FFN and Dense-PCA. The differences between FFN and Dense-PCA are fairly small; while Dense-PCA has a better mean value for predicting osmolarity, the difference is small compared to the MAD spreads. The Dense-PCAX predictions are more centrally concentrated but substantially less accurate in mean for both quantities, even relative to the spreads. 

\begin{figure}
    \centering
    \includegraphics[width=\textwidth]{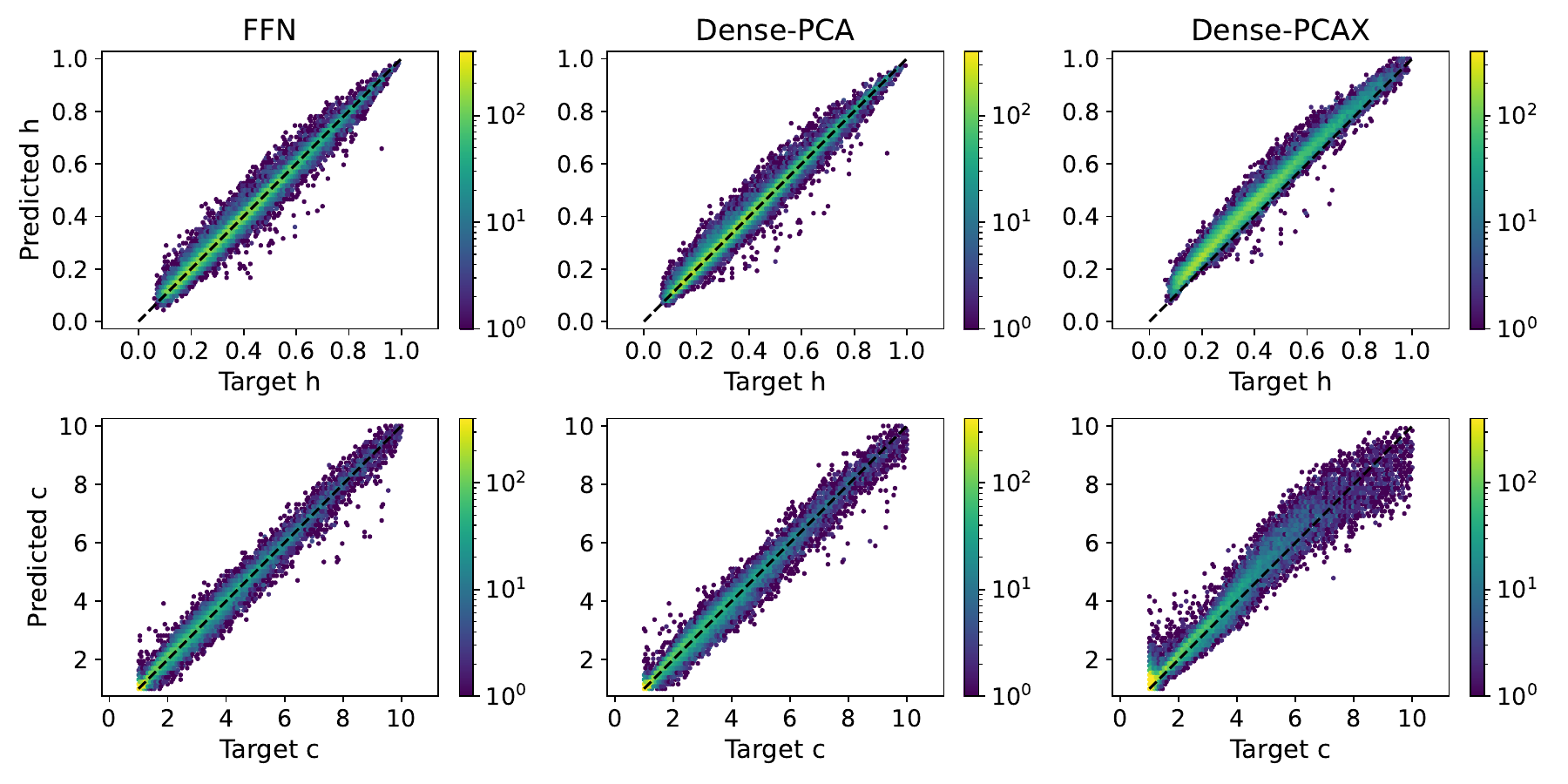}
    \caption{Comparisons of true and predicted final values of thickness (first row) and osmolarity (second row) for ODE-trained learners on all the synthetic ODE testing data. The color map shows the number of predictions in each bin.}
    \label{fig:FFNerr_ode_cl}
\end{figure}

\autoref{fig:FFNerr_ode_cl} compares all the individual ODE-trained predictions and true values of the final-time thickness $h(1)$ and osmolarity $c(1)$. For the prediction of $h(1)$, all three learners perform well, with Dense-PCAX exhibiting the tightest clustering about the diagonal, followed closely by FFN, while Dense-PCA shows slightly greater scatter. For the prediction of $c(1)$, FFN and Dense-PCA remain accurate across most of the range, whereas the poorer performance of Dense-PCAX is primarily due to underestimation when $c(1) > 8$.

\begin{figure}
  \centering
  \includegraphics[width=0.9\textwidth]{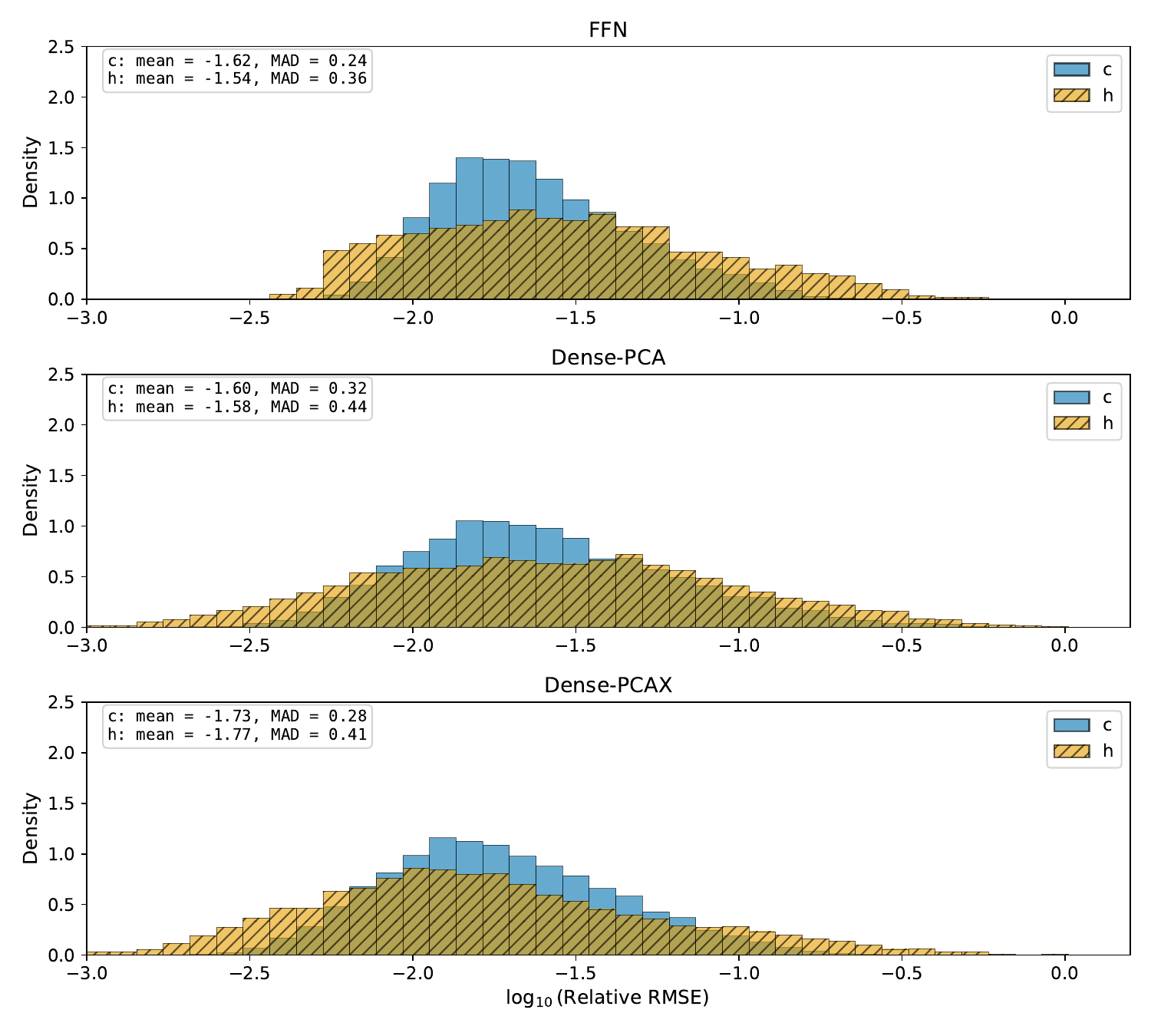}
    \caption{Relative RMSE for synthetic testing of ML learners trained by PDE data. The mean and mean absolute deviation (MAD) is shown for each distribution. Top to bottom: FFN, Dense-PCA, Dense-PCAX.}
%\caption{Relative RMSE errors in experimental data of three models. Top to Bottom: FFN model,PCA-dense model, and enhanced PCA model.}
  \label{fig:PDE_rmse_syn}
\end{figure}

\begin{figure}
    \centering
    \includegraphics[width=\textwidth]{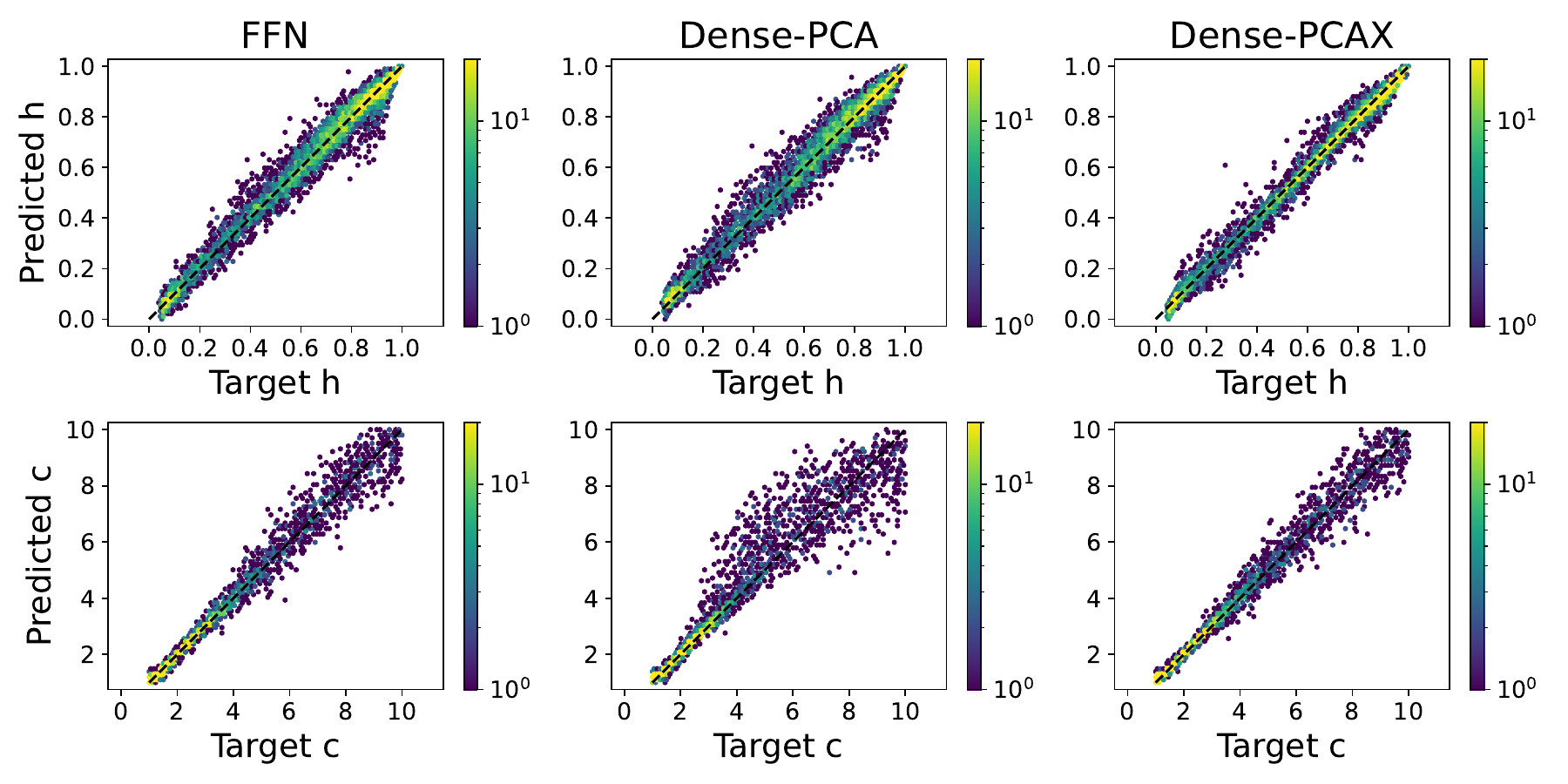}
    \caption{Comparisons of true and predicted final values of thickness (first row) and osmolarity (second row) for PDE-trained learners on all the synthetic PDE testing data. The color map shows the number of predictions in each bin.}
    \label{fig:PDE_final_c_syn}
\end{figure}

\autoref{fig:PDE_rmse_syn} shows the rRMSE values of osmolarity and thickness of the PDE-trained learners on the PDE test data. As with the ODE learners, the FFN and Dense-PCA distributions are similar to each other. Here, however, the Dense-PCAX learners do a bit better than the others on average, though by less than the MAD spread. \autoref{fig:PDE_final_c_syn} shows analogous results for the synthetic PDE test set. Again, all three learners accurately predict the final thickness, although Dense-PCA is slightly more scattered than FFN and Dense-PCAX. For the final osmolarity, Dense-PCA shows the greatest variance, while FFN and Dense-PCAX remain more tightly clustered around the diagonal. The prediction task is most difficult for all three learners when $c(1)$ is large.

\subsection{Experimental data}

The experimental data used in this study were collected at Indiana University as part of a study approved by the Biomedical Institutional Review Board \citep{awisi-gyauChangesCornealDetection2019a}. Fluorescein imaging data were obtained from 25 participants. Prior to enrollment, all subjects completed a screening procedure and subsequently attended two study visits, during which ten imaging trials were performed at each visit. Individuals who wore contact lenses or had received a clinical diagnosis of dry eye disease were excluded from participation.

Each trial consisted of a sequence of ocular surface images acquired at intervals of either $0.2$ or $0.25$ seconds, depending on the camera frame rate. At the start of each trial, a 2 microliter drop of 2\% sodium fluorescein solution was instilled into the eye \citep{carlson2016clinical}. The tear film was illuminated using a cobalt blue excitation filter (494 nm), and the emitted fluorescence was captured through a Wratten no. 8 yellow barrier filter positioned along the imaging axis. Under this illumination, fluorescein within the tear film emitted green light with a wavelength of approximately 521 nm, which was recorded by the imaging system \citep{book}.

To begin the recording, subjects were instructed to blink three times in order to distribute the fluorescein uniformly across the tear film. During these initial blinks, the illumination intensity was kept low, and a custom MATLAB algorithm \cite{wuEffectsIncreasingOcular2015} was used to estimate the initial fluorescein concentration, which was assumed to be spatially uniform over the corneal surface. Following the third blink, the illumination was increased to a predetermined high-intensity setting, and subjects were asked to keep their eyes open for as long as possible. Data collection continued until the next blink, which marked the end of the trial.

\citet{driscollFittingODEModels2023} fit 467 intensity time series of presumed TBU to the ODE model using traditional inverse problem optimization, thereby producing solutions $h$ and $c$ consistent with the observed intensity. Our ML models are meant to accomplish the same predictions without reference to the mathematical models once trained. Note that only cases for which a good ODE fit was found were reported in that work, introducing a form of bias, and that the possible existence of fits to $I$ nearly as good with different predictions of $h$ and $c$ is unknown. Moreover, while the ODE-trained predictors have learned the dynamics from the same mathematical models as in~\citet{driscollFittingODEModels2023}, the PDE-trained predictors presumably have a more physically realistic picture of the dynamics. Thus, while we continue to report "errors" in the following results with respect to~\citet{driscollFittingODEModels2023} as reference values, they might better be described as "discrepancies" for the PDE predictions.

\begin{figure}
    \centering
    \includegraphics[width=1.0\textwidth]{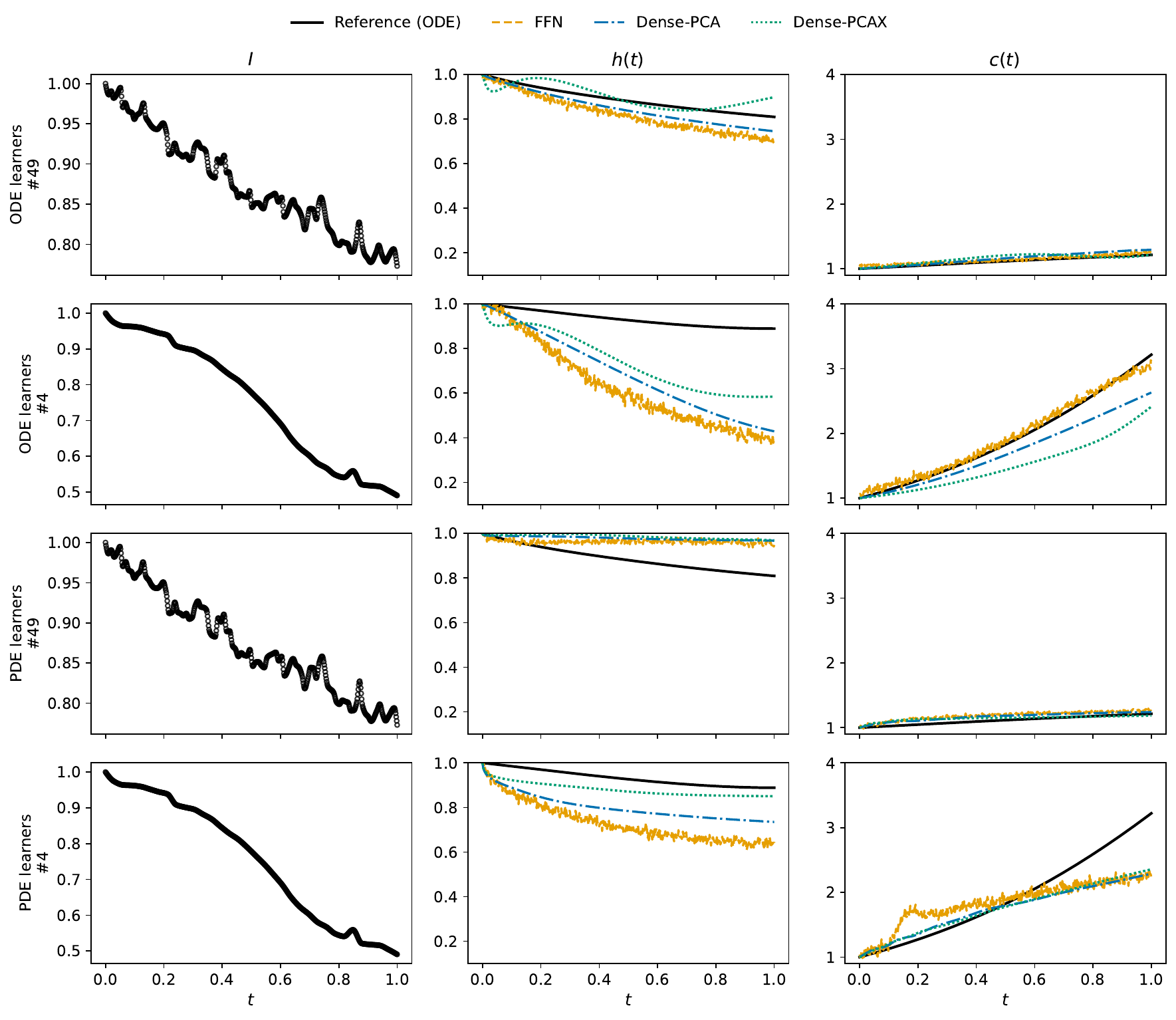}
    \caption{Example predictions on two time series from the experimental data by the ODE (top two rows) and PDE (bottom two rows) learners.}
    \label{fig:example-predictions-experimental}
\end{figure}

\autoref{fig:example-predictions-experimental} shows example predictions of all the learners for two of the experimental intensity inputs. Compared to predictions on synthetic data as seen in \autoref{fig:example-predictions-synthetic-ode} and \autoref{fig:example-predictions-synthetic-pde}, the predictions show more variance and less accuracy, particularly for thickness prediction, though conclusions cannot be drawn from only two examples. 

\begin{figure}
    \centering
    \includegraphics[width=0.8\textwidth]{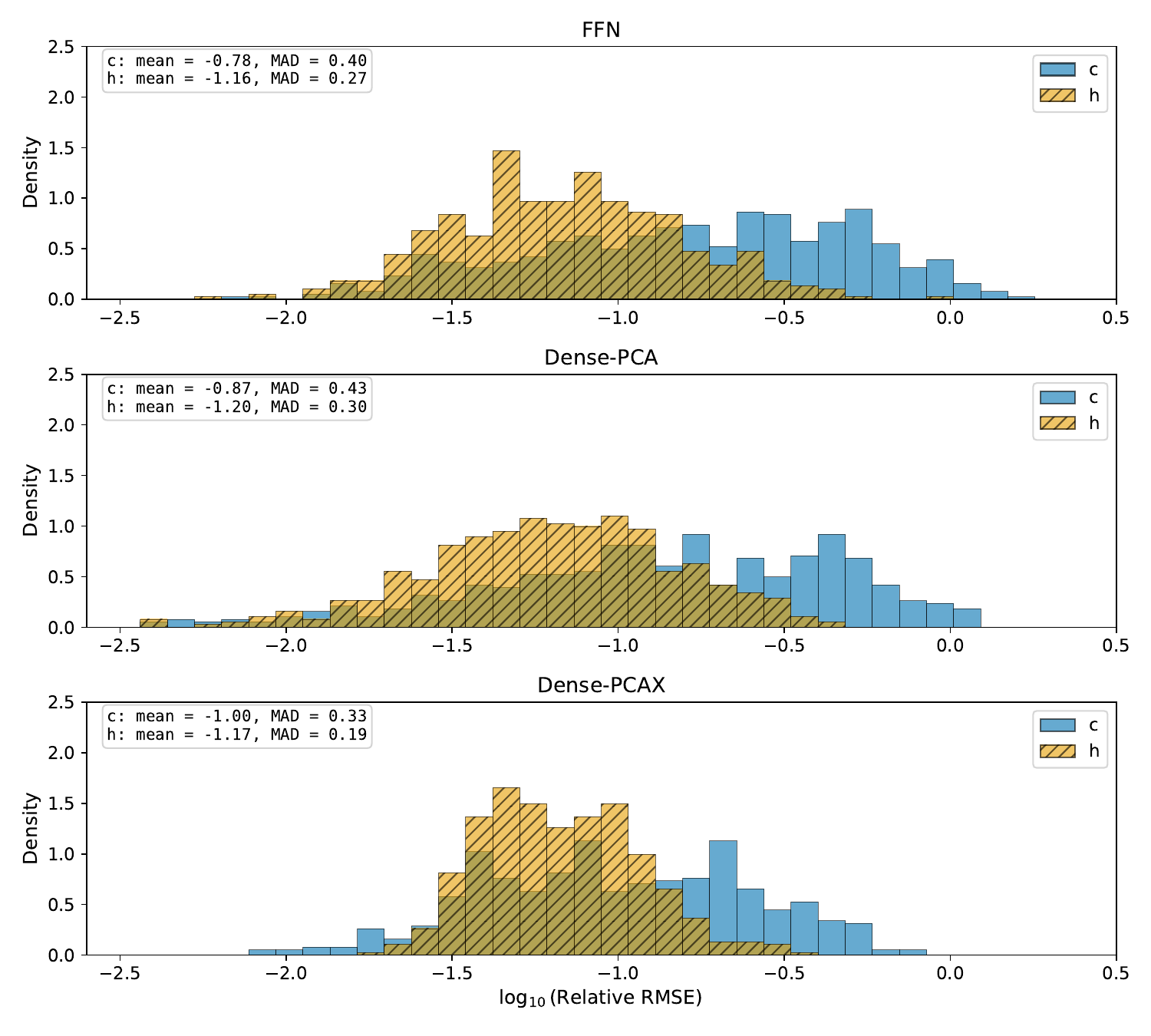}
    \caption{Relative RMSE for ODE-trained learners on the experimental data. }
    \label{fig:ode_exp_err}
\end{figure}

\begin{figure}
    \centering
    \includegraphics[width=0.8\textwidth]{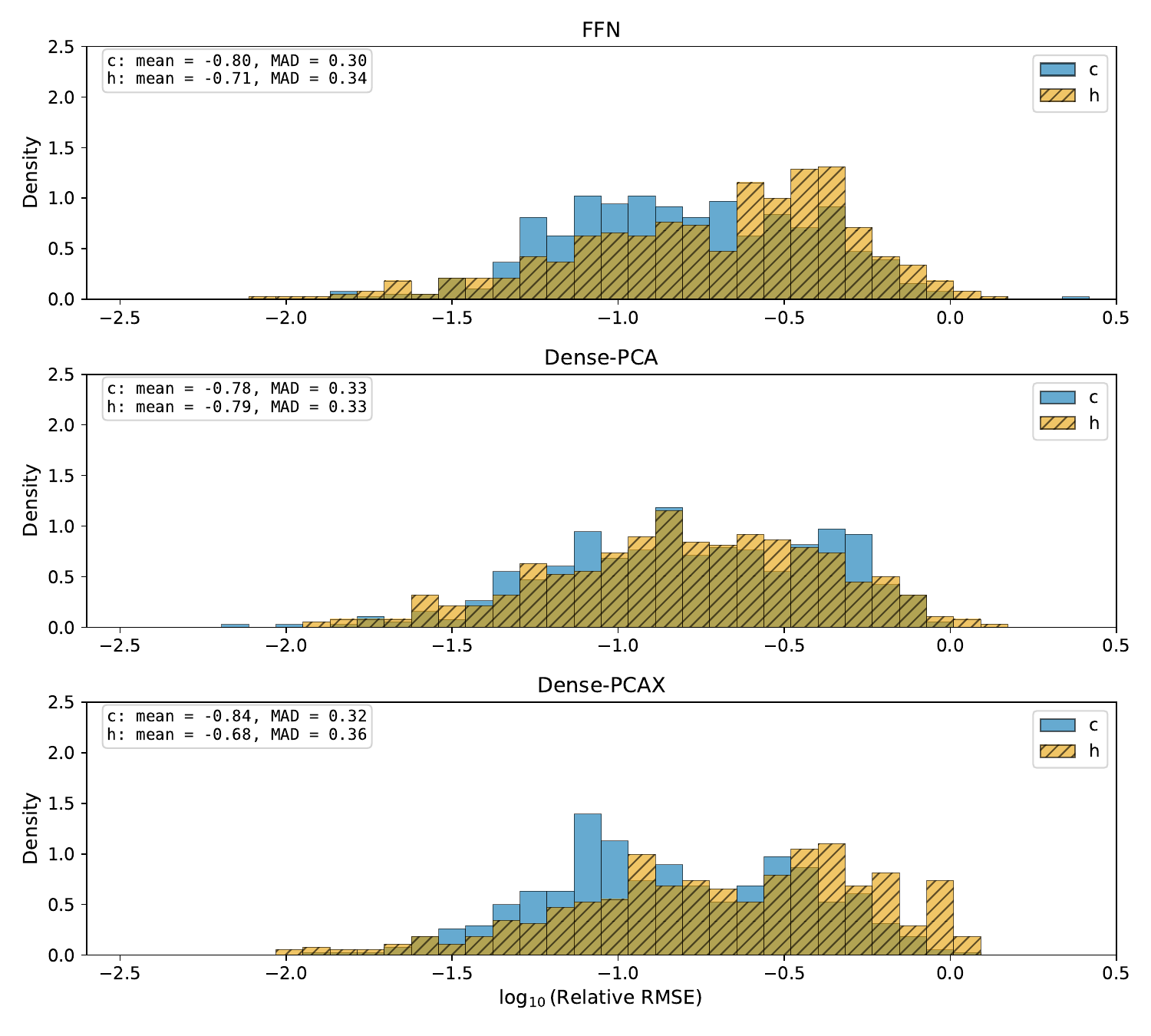}
    \caption{Relative RMSE for PDE-trained learners on the experimental data. }
    \label{fig:pde_exp_err}
\end{figure}

\autoref{fig:ode_exp_err} shows the rRMSE distributions of predictions by the ODE-trained learners on the experimental data. For these predictions, the Dense-PCAX is a little better than the others on average for osmolarity, but by an amount less than the MAD. The thickness distribution has less spread but is no better on average.  \autoref{fig:pde_exp_err} shows the distributions for PDE-trained learners. While their shapes vary, the quantitative aggregate differences are not large. 

\begin{figure}[tp]
    \centering
    \includegraphics[width=0.95\textwidth]{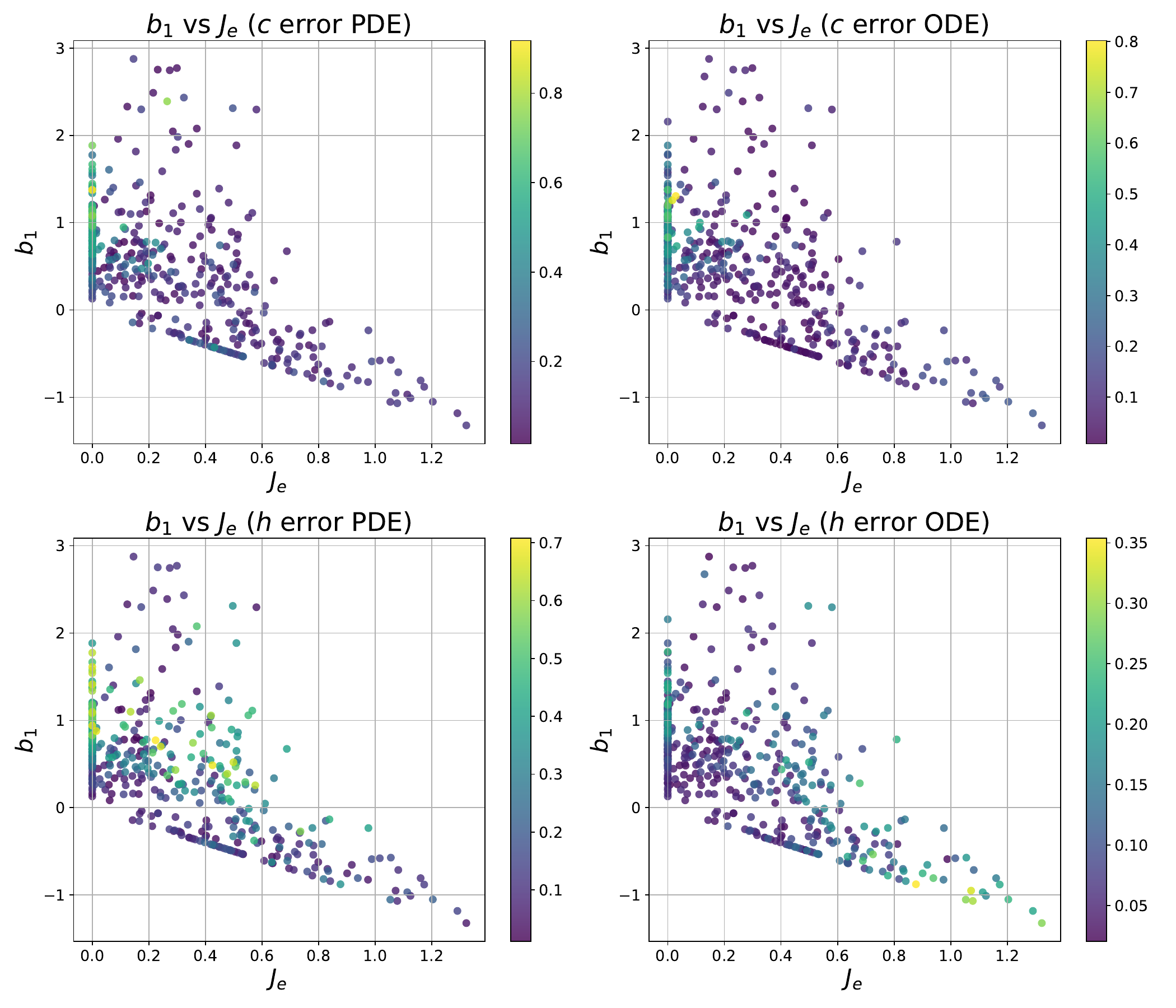}
    \caption{By rows: osmolarity error to thickness error. By columns: PDE Dense-PCAX learner to ODE Dense-PCAX learner. Color indicates the error comparing to the calculated values, shown as functions of the nondimensional ODE parameters $J_e$, $b_1$. }
    \label{fig:pde_exp_params}
\end{figure}

\autoref{fig:pde_exp_params} shows the dependence of rRMSE on the ODE fit parameters $b_1$ and $J_e$, which are proxies for shear flow rate and evaporation rate. The ODE-trained predictors show the greatest discrepancies in osmolarity for large $b_1>0$ and small $J_e$---that is, flow-dominated thinning---while the greatest discrepancies in thickness are for large $J_e$ and $b_1<0$, which indicates strong evaporation with resultant inward flow. For the PDE-trained predictor, osmolarity shows the same pattern as for the ODE-trained learner, while thickness discrepancy is largest for moderate outward flow over a wide range of evaporation rates.

\begin{figure}[tp]
    \centering
    \includegraphics[width=0.8\textwidth]{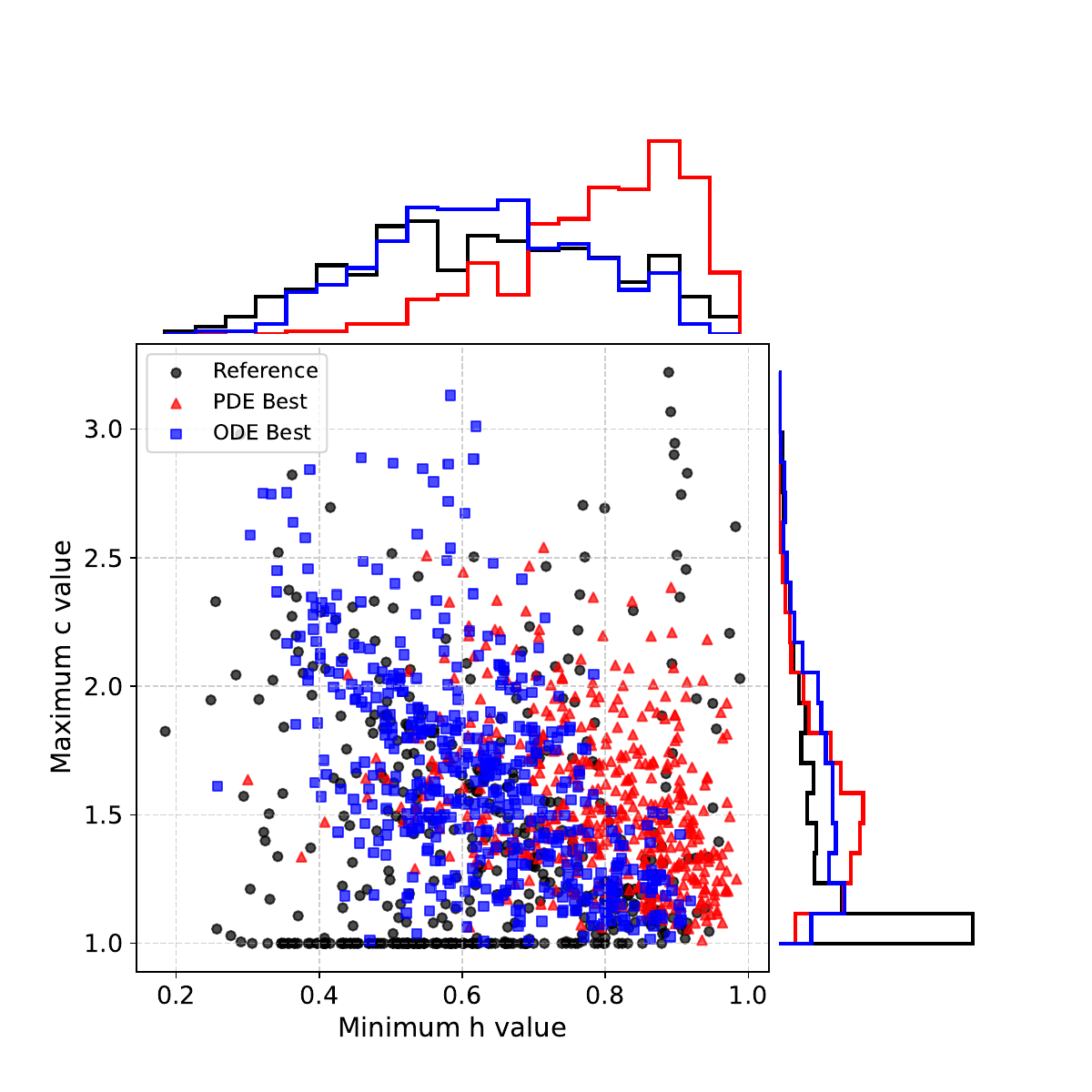}
    \caption{Predictions of maximum osmolarity and minimum thickness for~\citet{driscollFittingODEModels2023} (black circles) and two ML learners (red triangles for PDE training, blue squares for ODE training). The top and right histograms show the marginal distributions of thickness and osmolarity, respectively.}
    \label{fig:c_h_scatter_compare}
\end{figure}

Finally, \autoref{fig:c_h_scatter_compare} shows the predictions of minimum thickness and maximum osmolarity for the reference data together with the predictions from the best-performing learner for each test instance, selected separately from the ODE-trained and PDE-trained models. Among the ODE-trained learners, Dense-PCA was selected most frequently. Among the PDE-trained learners, Dense-PCAX was selected most frequently for predicting maximum osmolarity, whereas Dense-PCA was selected most frequently for predicting minimum thickness. The distribution of minimum thickness is similar for the ODE predictor and reference, but the reference predicts far more cases of maximum osmolarity staying close to 1. The PDE predictor is much more likely to predict minimum thickness greater than $0.75$ and a moderately larger maximum osmolarity on average.  

\section{Discussion}
\label{sec:discussion}

For this study, we used operator learning for the map between time series of the FL intensity and the film thickness and osmolarity. These outputs have clear physical significance, while the parameters of a particular mathematical model have a more indirect interpretation. However, one could instead train for the functional mapping to model parameter values and then simulate the model to obtain the outputs, as in a traditional inverse method, if the parameter values themselves are of interest.

In a broad sense, the primary task of learning the operator is to distinguish between thinning due to evaporation and thinning due to flow. Evaporative thinning tends to be slow---modeled mathematically as linear decrease---while thinning due to shear flow is exponential in nature, but the rates need not be clearly distinct in every case. Distinguishing between these effects is complicated by the interaction between $h$ and $c$ in producing the FL intensity $I$. 

Multiple operator learning methods were trained on datasets comprising simulations of ODE models and of PDE models. Performance, as measured on synthetic test data, did not depend strongly on the particular learning method or architecture used: 1--2 digits of rRMSE accuracy were found in most cases. Given the weak dependence on the learning method, we speculate that the accuracy of the predictions is limited mainly by the inherent identifiability of the dependent quantities from FL data presented in this form. 

Compared to the inverse problem predictions of~\citet{driscollFittingODEModels2023} on experimental data, the ML predictions based on the same ODE model are much less likely to predict minor change in osmolarity simultaneously with substantial thinning, and large osmolarity simultaneously with little thinning. Differences in osmolarity predictions appear most often in what~\citet{driscollFittingODEModels2023} describes as flow-dominated, zero-evaporation instances, while differences in thickness predictions also occur for strongly evaporative instances. One possibility is that the standard ODE inverse problem can be more easily induced to produce solution types that play a small role in the full solution set. Since the underlying physical realities are unknown, we cannot determine whether such results represent appropriate inverse problem responses or a manifestation of overfitting in the standard inverse problem.

On experimental data, the PDE-trained learners predict much less thinning on average than either of the ODE-based predictors do over a wide range of evaporation rates. This suggests a fundamental difference between the ODE and PDE models of the phenomenon. While the PDE predictions have even fewer cases than the ODE learners of osmolarity remaining very close to 1, they do make more predictions between 1.25 and 1.75 and thus are distributed more like the inverse problem reference results. It might be informative to compare PDE operator learning to traditional inverse-problem methods for the same PDE on the data, but that project is beyond the scope of this study.

It is plausible that PDE-trained learners would perform better if the training data and observations included information from more than a single location on the cornea. This would substantially increase the size of inputs and outputs to the learners, increasing training times. We also note that the learners in this work have no constraints from physics or dynamics other than what is implied by the training data. Methods that can incorporate such information, such as physics-informed neural nets~\citep{CuomoScientificMachine2022,LiPhysicsInformedNeural2023}, neural ODEs~\citep{ChenNeuralOrdinary2019,KidgerNeuralDifferential2022}, and universal differential equations~\citep{RackauckasUniversalDifferential2020} might be expected to perform better. Finally, models and experimental observations that included additional variables might show better results. In particular, thermal imaging is experimentally practical and could give a strong signal through evaporative cooling that is independent of flow.

\vspace{1em} \noindent \textbf{Conflict of Interest}. The authors declare no competing interests.

\bibliographystyle{abbrvnat}
\bibliography{refs.bib, sciml_cite.bib}
\end{document}